\documentclass[11pt]{article}
\usepackage{amsfonts}
\usepackage{amsfonts}
\usepackage{amsfonts}
\usepackage{mathrsfs}
\usepackage{bm}
\usepackage{cite}
\usepackage[colorlinks,linkcolor=blue,citecolor=blue]{hyperref}
\usepackage{amssymb,amsmath} 
 \allowdisplaybreaks

\catcode`!=11
\let\!int\int \def\int{\displaystyle\!int}
\let\!lim\lim \def\lim{\displaystyle\!lim}
\let\!sum\sum \def\sum{\displaystyle\!sum}
\let\!sup\sup \def\sup{\displaystyle\!sup}
\let\!inf\inf \def\inf{\displaystyle\!inf}
\let\!cap\cap \def\cap{\displaystyle\!cap}
\let\!max\max \def\max{\displaystyle\!max}
\let\!min\min \def\min{\displaystyle\!min}
\let\!frac\frac \def\frac{\displaystyle\!frac}
\catcode`!=12

\let\oldsection\section
\renewcommand\section{\setcounter{equation}{0}\oldsection}

\allowdisplaybreaks
\def\pf{\it{Proof.}\rm\quad}

\def\N{\mathbb{N}}

\newtheorem{thm}{Theorem}[section]
\newtheorem{lem}[thm]{Lemma}
\newtheorem{cor}[thm]{Corollary}

\begin{document}
\title {\bf Some evaluation of cubic Euler sums}
\author{
{Ce Xu\thanks{Corresponding author. Email: 15959259051@163.com}}\\[1mm]
\small School of Mathematical Sciences, Xiamen University\\
\small Xiamen
361005, P.R. China}

\date{}
\maketitle \noindent{\bf Abstract} P. Flajolet and B. Salvy \cite{FS1998} prove the famous theorem that a nonlinear Euler sum $S_{i_1i_2\cdots i_r,q}$ reduces to a combination of sums of lower orders whenever the weight $i_1+i_2+\cdots+i_r+q$ and the order $r$ are of the same parity. In this article, we develop an approach to evaluate the cubic sums $S_{1^2m,p}$ and $S_{1l_1l_2,l_3}$. By using the approach, we establish some relations involving cubic, quadratic and linear Euler sums. Specially, we prove the cubic sums $S_{1^2m,m}$ and $S_{1(2l+1)^2,2l+1}$ are reducible to zeta values, quadratic and linear sums. Moreover, we prove that the two combined sums involving multiple zeta values of depth four
\[\sum\limits_{\left\{ {i,j} \right\} \in \left\{ {1,2} \right\},i \ne j} {\zeta \left( {{m_i},{m_j},1,1} \right)}\quad {\rm and}\quad \sum\limits_{\left\{ {i,j,k} \right\} \in \left\{ {1,2,3} \right\},i \ne j \ne k} {\zeta \left( {{m_i},{m_j},{m_k},1} \right)} \]
can be expressed in terms of multiple zeta values of depth $\leq 3$, here $2\leq m_1,m_2,m_3\in \N$.
Finally, we evaluate the alternating cubic Euler sums ${S_{{{\bar 1}^3},2r + 1}}$ and show that it are reducible to alternating quadratic and linear Euler sums. The approach is based on Tornheim type series computations.
\\[2mm]
\noindent{\bf Keywords} Harmonic number; polylogarithm function; Euler sum; Tornheim type series; Riemann zeta function, multiple zeta value.
\\[2mm]
\noindent{\bf AMS Subject Classifications (2010):} 11M06; 11M32; 11M99; 33E20; 40A05.
\tableofcontents
\section{Introduction}
In response to a letter from Goldbach, Euler considered sums of the form (see Berndt \cite{B1985} for a discussion)
\begin{align}\label{1.1}
S_{p,q}:=\sum\limits_{n = 1}^\infty  {\frac{{H_n^{\left( p \right)}}}
{{{n^q}}}},
\end{align}
and was able to give explicit values for certain of these sums in terms of the Riemann zeta function, where $p,q$ are positive integers with $q \geq 2$, and $w:=p+q$ denotes the weight of linear sums $S_{p,q}$. These kind of sums are called the linear Euler sums (for short linear sums) today.
Here $H^{(p)}_n$ denotes the harmonic number which is defined by
\[
H^{(k)}_n:=\sum\limits_{j=1}^n\frac {1}{j^k}\quad {\rm and}\quad H^{(k)}_0:=0.
\]
When $k=1$, then $H_n:=H^{(1)}_n$, which is called the classical harmonic number.

In their famous paper \cite{FS1998}, Flajolet and Salvy introduced the following generalized series
\begin{align}\label{1.2}
{S_{{\bf S},q}} := \sum\limits_{n = 1}^\infty  {\frac{{H_n^{\left( {{s_1}} \right)}H_n^{\left( {{s_2}} \right)} \cdots H_n^{\left( {{s_r}} \right)}}}
{{{n^q}}}},
\end{align}
which is called the generalized (nonlinear) Euler sums. Here ${\bf S}:=(s_1,s_2,\ldots,s_r)\ (r,s_i\in \N, i=1,2,\ldots,r)$ and $q\geq 2$. The quantity $w:={s _1} +  \cdots  + {s _r} + q$ is called the weight and the quantity $r$ is called the degree. As usual, repeated summands in partitions are indicated by powers, so that for instance
\[{S_{{1^2}{2^3}4,q}} = {S_{112224,q}} = \sum\limits_{n = 1}^\infty  {\frac{{H_n^2[H^{(2)} _n]^3{H^{(4)} _n}}}{{{n^q}}}}. \]

It has been discovered in the course of the years that many Euler sums admit expressions involving finitely the ``zeta values", that is to say of the Riemann zeta function \cite{A2000},
$$\zeta(s):=\sum\limits_{n = 1}^\infty {\frac {1}{n^{s}}},\Re(s)>1$$
with positive integer arguments. For example, the linear sums $S_{p,q}$ can be evaluated in terms of zeta values in the following cases: $p=1,p=q,p+q$ odd and $p+q=6$ with $q\geq 2$\ (for more details, see \cite{BBG1994,BBG1995,FS1998}). In 1994, Bailey et al. \cite{BBG1994} proved that all Euler sums of the form $S_{1^p,q}$ for weights $p+q\in \{3,4,5,6,7,9\}$ are reducible to {\bf Q}-linear combinations of zeta values by using the experimental method. In \cite{X2017,Xu2017}, we proved that all Euler sums of weight $\leq 7$ are reducible to $\mathbb{Q}$-linear combinations of single zeta monomials. For weight 9, all Euler sums of the form ${S_{{s_1} \cdots {s_k},q}}$ with $q\in \{4,5,6,7\}$ are expressible polynomially in terms of zeta values. Very recently, Wang et al \cite{W2017} shown that all Euler sums of weight nine are reducible to zeta values.

However, there are also many nonlinear Euler sums which need not only zeta values but also linear sums. Namely, many nonlinear Euler sums are reducible to polynomials in zeta values and to linear sums (see \cite{BBG1994,BBG1995,FS1998,M2014,Xu2016,X2016}). For instance, in 1995, Borwein et al. \cite{BBG1995} showed that the quadratic sums $S_{ 1^2,q}$ can reduce to linear sums $S_{2,q}$ and polynomials in zeta values. In 1998, Flajolet and Salvy \cite{FS1998} used the contour integral representations and residue computation to show that the quadratic sums $S_{p_1p_2,q}$ are reducible to linear sums and zeta values when the weight $p_1 + p_2 + q$ is even and $p_1,p_2>1$. In \cite{Xu2017}, we proved that all Euler sums with weight eight are reducible to zeta values and linear sum $S_{2,6}$.

Hence, a good deal of work on Euler sums has focused on the problem of determining when `complicated' sums can be expressed in terms of `simpler' sums. Thus, researchers are interested in determining which sums can be expressed in terms of other sums of lesser degree. Besides the works referred to above, there are many other researches devoted to the Euler sums. For example, please see \cite{A1997,BBGP1996,BZB2008,CB1994,ELY2005,F2005,M2014} and references therein.

Multiple zeta values (also called Zagier sums) are the natural generalizations of classical Euler sums from double sums to more general $k$-fold Euler sums defined by \cite{BBBL1997,BBBL1999,Xu2017}
\begin{align}\label{1.3}
\zeta \left( \mathbf{S} \right) \equiv \zeta \left( {{s_1},{s_2}, \cdots ,{s_k}} \right) := \sum\limits_{{n_1} > {n_2} >  \cdots  > {n_k} \ge 1} {\frac{1}{{n_1^{{s_1}}n_2^{{s_2}} \cdots n_k^{{s_k}}}}},
\end{align}
where $s_1,s_2,\ldots,s_k$ are positive integers with $s_1>1$ for the sake of convergence. The numbers $w:={s_1} +  \cdots  + {s_k}$ and $k$ are called the weight and depth of $\zeta \left( {{s_1},{s_2}, \cdots ,{s_k}} \right)$  respectively.
For convenience, we let $\{a\}_k$ be the $k$ repetitions of a such that
\[\zeta \left( {5,3,{{\left\{ 1 \right\}}_4},2} \right) = \zeta \left( {5,3,1,1,1,1,2} \right).\]
Many papers use the opposite convention, with the $n_i$'s ordered by $n_1<n_2<\cdots<n_k$ or $n_1\leq n_2\leq \cdots\leq n_k$, see \cite{H1992,DZ2012}.

From the definitions of Euler sums and multiple zeta values, we can find the relations
\begin{align}\label{1.4}
{S_{p,q}} = \zeta \left( {q,p} \right) + \zeta \left( {q + p} \right),
\end{align}
\begin{align}\label{1.5}
&{S_{{p_1}{p_2},q}} = \zeta \left( {q,{p_1},{p_2}} \right) + \zeta \left( {q,{p_2},{p_1}} \right) + \zeta \left( {q + {p_1},{p_2}} \right)\nonumber\\
&\quad\quad\quad\quad+ \zeta \left( {q + {p_2},{p_1}} \right) + \zeta \left( {q,{p_1} + {p_2}} \right) + \zeta \left( {q + {p_1} + {p_2}} \right).
\end{align}
Multiple zeta values were introduced and studied by Euler in the old days. The multiple zeta values have attracted
considerable interest in recent years. In the past two decades, many authors have studied multiple zeta values, and a number of relations among them have been found, see \cite{BBBL1997,BBBL1999,BG1996,EW2012,H1992,O2008,O2010,M2000,T2004,Xu2017,DZ1994,DZ2012,Z2015} and references therein. For example, some results on triple zeta values and quadruple zeta values were evaluated in \cite{BG1996,ELY2005,EW2012,O2008,O2010}. Tsumura \cite{T2004} proved in 2004 that the multiple zeta value  $\zeta \left( {{s_1},{s_2}, \cdots ,{s_k}} \right)$
can be expressed as a rational linear combination of products of multiple zeta values of lower depth provided that the sum of weight and depth, $w+k$, is odd. Jonathan M. Borwein, David M. Bradley and David J. Broadhurst \cite{BBBL1997} proved the following duality relation
\begin{align}\label{1.6}
\zeta \left( {{m_1} + 2,{{\left\{ 1 \right\}}_{{n_1}}}, \ldots ,{m_p} + 2,{{\left\{ 1 \right\}}_{{n_p}}}} \right) = \zeta \left( {{n_p} + 2,{{\left\{ 1 \right\}}_{{m_p}}}, \ldots ,{n_1} + 2,{{\left\{ 1 \right\}}_{{m_1}}}} \right).
\end{align}
Zagier \cite{DZ2012} proved that the multiple zeta values ${\zeta \left( {{{\left\{ 2 \right\}}_a},3,{{\left\{ 2 \right\}}_b}} \right)} $ are reducible to polynomials in zeta values, $a,b\in \N_0:=\{0,1,2,\ldots\}$, and gave explicit formulae.

In this paper, we develop an approach to evaluation of Euler sums. The approach is based on Tornheim type
series computations. The purpose of this paper is to establish some explicit formulas of Euler sums with weight even and degree three in terms of zeta values, linear and quadratic sums by using the method of Tornheim type series computations. Moreover, by using the approach, we can obtain some relations between Euler sums of degree three and multiple zeta values of depth four.

\section{Main Theorems and Corollaries}
The main theorems and corollaries of this paper can be stated as follows.
\begin{thm}\label{thm1} For positive integers $m,p\geq 2$, then the following identity holds:
\begin{align}\label{2.1}
 &{S_{{1^2}m,p}} + {S_{{1^2}p,m}} + {S_{2m,p}} + {S_{2p,m}}\nonumber\\
  &= 2{\left( { - 1} \right)^{p - 1}}\sum\limits_{i = 1}^{p - 2} {{{\left( { - 1} \right)}^{i - 1}}{S_{m,i + 1}} \cdot {S_{1,p - i}}}  + 2{\left( { - 1} \right)^{m - 1}}\sum\limits_{i = 1}^{m - 2} {{{\left( { - 1} \right)}^{i - 1}}{S_{p,i + 1}} \cdot {S_{1,m - i}}} \nonumber \\
  &\quad- 2\zeta \left( 2 \right)\left( {\zeta \left( m \right)\zeta \left( p \right) + \zeta \left( {m + p} \right)} \right) + 2{S_{1m,p + 1}} + 2{S_{1p,m + 1}}\nonumber \\
  &\quad+ 2{\left( { - 1} \right)^{m + p}}{S_{{1^2},m + p}} - 2{\left( { - 1} \right)^{m + p}}{S_{1,p}} \cdot {S_{1,m}} \nonumber\\
  &\quad- 2{\left( { - 1} \right)^{p - 1}}\zeta \left( {m + 1} \right){S_{1,p}} - 2{\left( { - 1} \right)^{m - 1}}\zeta \left( {p + 1} \right){S_{1,m}}\nonumber \\
 &\quad - 2{\left( { - 1} \right)^{p - 1}}\sum\limits_{j = 1}^{m - 1} {{{\left( { - 1} \right)}^{j - 1}}\zeta \left( {m + 1 - j} \right)\left( {{S_{1j,p}} - {S_{1,p + j}}} \right)} \nonumber \\
 &\quad - 2{\left( { - 1} \right)^{m - 1}}\sum\limits_{j = 1}^{p - 1} {{{\left( { - 1} \right)}^{j - 1}}\zeta \left( {p + 1 - j} \right)\left( {{S_{1j,m}} - {S_{1,m + j}}} \right)} .
\end{align}
\end{thm}
\begin{thm}\label{thm2} If $l_1+l_2+l_3$ is odd, and $l_1,l_2,l_2\geq 2$ are positive integers, then the cubic combination $S_{1l_1l_2,l_3}+S_{1l_1l_3,l_2}+S_{1l_2l_3,l_1}$ is expressible in terms of quadratic, linear sums and zeta values. We have
\begin{align}
 &{S_{1{l_1}{l_2},{l_3}}} + {S_{1{l_1}{l_3},{l_2}}} + {S_{1{l_2}{l_3},{l_1}}} \nonumber\\
  & = {\left( { - 1} \right)^{{l_3}}}\sum\limits_{i = 1}^{{l_2} - 2} {{{\left( { - 1} \right)}^{i - 1}}{S_{{l_3},i + 1}} \cdot {S_{{l_1},{l_2} - i}}}  + {S_{1,{l_1}}}\left( {\zeta \left( {{l_2}} \right)\zeta \left( {{l_3}} \right) + \zeta \left( {{l_2} + {l_3}} \right)} \right)  \nonumber\\
  &\quad+ {\left( { - 1} \right)^{{l_2}}}\sum\limits_{i = 1}^{{l_1} - 2} {{{\left( { - 1} \right)}^{i - 1}}{S_{{l_2},i + 1}} \cdot {S_{{l_3},{l_1} - i}}}  + {S_{1,{l_2}}}\left( {\zeta \left( {{l_1}} \right)\zeta \left( {{l_3}} \right) + \zeta \left( {{l_1} + {l_3}} \right)} \right)  \nonumber\\
  &\quad+ {\left( { - 1} \right)^{{l_1}}}\sum\limits_{i = 1}^{{l_3} - 2} {{{\left( { - 1} \right)}^{i - 1}}{S_{{l_1},i + 1}} \cdot {S_{{l_2},{l_3} - i}}}  + {S_{1,{l_3}}}\left( {\zeta \left( {{l_1}} \right)\zeta \left( {{l_2}} \right) + \zeta \left( {{l_1} + {l_2}} \right)} \right)  \nonumber\\
  &\quad+ {\left( { - 1} \right)^{{l_1} - 1}}\sum\limits_{j = 1}^{{l_1} - 1} {{{\left( { - 1} \right)}^{j - 1}}\zeta \left( {{l_1} + 1 - j} \right)\left( {{S_{{l_3}j,{l_2}}} + {S_{{l_2}j,{l_3}}} - {S_{{l_3},{l_2} + j}} - {S_{{l_2},{l_3} + j}}} \right)}  \nonumber\\
  &\quad + {\left( { - 1} \right)^{{l_2} - 1}}\sum\limits_{j = 1}^{{l_2} - 1} {{{\left( { - 1} \right)}^{j - 1}}\zeta \left( {{l_2} + 1 - j} \right)\left( {{S_{{l_3}j,{l_1}}} + {S_{{l_1}j,{l_3}}} - {S_{{l_3},{l_1} + j}} - {S_{{l_1},{l_3} + j}}} \right)} \nonumber\\
  &\quad+ {\left( { - 1} \right)^{{l_3} - 1}}\sum\limits_{j = 1}^{{l_3} - 1} {{{\left( { - 1} \right)}^{j - 1}}\zeta \left( {{l_3} + 1 - j} \right)\left( {{S_{{l_1}j,{l_2}}} + {S_{{l_2}j,{l_1}}} - {S_{{l_1},{l_2} + j}} - {S_{{l_2},{l_1} + j}}} \right)} \nonumber\\
  &\quad + {\left( { - 1} \right)^{{l_1} - 1}}\zeta \left( {{l_1} + 1} \right)\left( {\zeta \left( {{l_2}} \right)\zeta \left( {{l_3}} \right) + \zeta \left( {{l_2} + {l_3}} \right)} \right) \nonumber\\
  &\quad+ {\left( { - 1} \right)^{{l_2} - 1}}\zeta \left( {{l_2} + 1} \right)\left( {\zeta \left( {{l_1}} \right)\zeta \left( {{l_3}} \right) + \zeta \left( {{l_1} + {l_3}} \right)} \right) \nonumber\\
  &\quad+ {\left( { - 1} \right)^{{l_3} - 1}}\zeta \left( {{l_3} + 1} \right)\left( {\zeta \left( {{l_1}} \right)\zeta \left( {{l_2}} \right) + \zeta \left( {{l_1} + {l_2}} \right)} \right)\nonumber\\
  &\quad - {S_{1\left( {{l_1} + {l_2}} \right),{l_3}}} - {S_{1\left( {{l_1} + {l_3}} \right),{l_2}}} - {S_{1\left( {{l_2} + {l_3}} \right),{l_1}}}.
\end{align}
\end{thm}
\begin{thm}\label{thm3} For positive integers $m_1,m_2\geq 1$ and $m_3>1$, then the following relation holds:
\begin{align}\label{2.3}
{S_{1{m_1}{m_2},{m_3}}} = &{\left( { - 1} \right)^{{m_3} - 1}}\sum\limits_{i = 1}^{{m_3} - 2} {{{\left( { - 1} \right)}^{i - 1}}\zeta \left( {{m_3} - i} \right){S_{{m_1}{m_2},i + 1}}}\nonumber \\
& + {\left( { - 1} \right)^{{m_3}}}\left( {{S_{{m_1},{m_2} + 1}} + {S_{{m_2},{m_1} + 1}} - \zeta \left( {{m_1} + {m_2} + 1} \right)} \right)\zeta \left( {{m_3}} \right)\nonumber \\
& + {\left( { - 1} \right)^{{m_3}}}\left\{ {\sum\limits_{i = 1}^{{m_2} - 1} {{{\left( { - 1} \right)}^{i - 1}}{S_{{m_1},{m_2} + 1 - i}}\zeta \left( {{m_3},i} \right)}  + \sum\limits_{i = 1}^{{m_1} - 1} {{{\left( { - 1} \right)}^{i - 1}}{S_{{m_2},{m_1} + 1 - i}}\zeta \left( {{m_3},i} \right)} } \right\}\nonumber \\
& + {\left( { - 1} \right)^{{m_3}}}\left\{ {{{\left( { - 1} \right)}^{{m_2} - 1}}\zeta \left( {{m_1} + 1} \right)\zeta \left( {{m_3},{m_2}} \right) + {{\left( { - 1} \right)}^{{m_1} - 1}}\zeta \left( {{m_2} + 1} \right)\zeta \left( {{m_3},{m_1}} \right)} \right\}\nonumber \\
& + {\left( { - 1} \right)^{{m_3}}}\left\{ \begin{array}{l}
 {\left( { - 1} \right)^{{m_2} - 1}}\sum\limits_{i = 1}^{{m_1} - 1} {{{\left( { - 1} \right)}^{i - 1}}\zeta \left( {{m_1} + 1 - i} \right)\zeta \left( {{m_3},{m_2},i} \right)}  \\
  + {\left( { - 1} \right)^{{m_1} - 1}}\sum\limits_{i = 1}^{{m_2} - 1} {{{\left( { - 1} \right)}^{i - 1}}\zeta \left( {{m_2} + 1-i} \right)\zeta \left( {{m_3},{m_1},i} \right)}  \\
 \end{array} \right\}\nonumber \\
& +{\left( { - 1} \right)^{{m_1} + {m_2} + {m_3}}}\left\{ \begin{array}{l}
 \zeta \left( {{m_3},{m_2},{m_1},1} \right) + \zeta \left( {{m_3},{m_1},{m_2},1} \right) \\
  \quad+ \zeta \left( {{m_3},{m_2},{m_1} + 1} \right) + \zeta \left( {{m_3},{m_1},{m_2} + 1} \right) \\
  \quad+ \zeta \left( {{m_3},{m_1} + {m_2},1} \right) + \zeta \left( {{m_3},{m_1} + {m_2} + 1} \right) \\
 \end{array} \right\}\nonumber \\
 & - {\left( { - 1} \right)^{{m_3}}}\sum\limits_{i = 1}^{{m_1} + {m_2} - 1} {{{\left( { - 1} \right)}^{i - 1}}\zeta \left( {{m_1} + {m_2} + 1 - i} \right)\zeta \left( {{m_3},i} \right)} .
\end{align}
\end{thm}
Letting $p=m$ in Theorem \ref{thm1} and $l_1=l_2=l_3=2l+1$ in Theorem \ref{thm2}, we can get the following two corollaries.
\begin{cor}\label{cor3} For positive integer $m>1$, then the cubic sums
\[{S_{{1^2}m,,m}}= \sum\limits_{n = 1}^\infty  {\frac{{H_n^2H_n^{\left( m \right)}}}{{{n^m}}}} \]
are reducible to quadratic and linear sums. We have
\begin{align}
 {S_{{1^2}m,,m}} &= 2{\left( { - 1} \right)^{m - 1}}\sum\limits_{i = 1}^{m - 2} {{{\left( { - 1} \right)}^{i - 1}}{S_{m,i + 1}} \cdot {S_{1,m - i}}}  - \zeta \left( 2 \right)\left( {{\zeta ^2}\left( m \right) + \zeta \left( {2m} \right)} \right)\nonumber\\
  &\quad - 2{\left( { - 1} \right)^{m - 1}}\sum\limits_{j = 1}^{m - 1} {{{\left( { - 1} \right)}^{j - 1}}\zeta \left( {m + 1 - j} \right)\left( {{S_{1j,m}} - {S_{1,m + j}}} \right)}  \nonumber\\
  &\quad + 2{S_{1m,m + 1}} + {S_{{1^2},2m}} - S_{1,m}^2 - 2{\left( { - 1} \right)^{m - 1}}\zeta \left( {m + 1} \right){S_{1,m}} - {S_{2m,m}}.
\end{align}
\end{cor}
\begin{cor}\label{cor4} For positive integer $l$, then the cubic sums
\[{S_{1{{\left( {2l + 1} \right)}^2},\left( {2l + 1} \right)}} = \sum\limits_{n = 1}^\infty  {\frac{{{H_n}\left(H^{(2l+1)} _n\right)^2}}{{{n^{2l + 1}}}}} \]
are reducible to quadratic and linear sums. We have
\begin{align}
 {S_{1{{\left( {2l + 1} \right)}^2},\left( {2l + 1} \right)}} &= \zeta \left( {2l + 2} \right)\left( {{\zeta ^2}\left( {2l + 1} \right) + \zeta \left( {4l + 2} \right)} \right) + \left( {\zeta \left( {4l + 2} \right) + {\zeta ^2}\left( {2l + 1} \right)} \right){S_{1,2l + 1}} \nonumber\\
  &\quad - {\left( { - 1} \right)^l}S_{2l + 1,l + 1}^2 - {S_{1\left( {4l + 2} \right),2l + 1}} - 2\sum\limits_{i = 1}^l {{{\left( { - 1} \right)}^{i - 1}}{S_{2l + 1,i + 1}} \cdot {S_{2l + 1,2l + 1 - i}}}  \nonumber\\
  &\quad+ 2\sum\limits_{j = 1}^{2l} {{{\left( { - 1} \right)}^{j - 1}}\zeta \left( {2l + 2 - j} \right)} \left( {{S_{j\left( {2l + 1} \right),2l + 1}} - {S_{2l + 1,2l + j + 1}}} \right).
\end{align}
\end{cor}
\section{Proofs of Theorem \ref{thm1}, \ref{thm2} and \ref{thm3}}
In this section, we will use the Tornheim type series $T_1\left( {l,m;p} \right)$ and $T_2\left( {l,m;p} \right)$ to prove the Theorem \ref{thm1}, \ref{thm2} and Theorem \ref{thm3}, respectively. Here the Tornheim type series $T_1\left( {l,m;p} \right)$ and $T_2\left( {l,m;p} \right)$ are defined by
\begin{align}\label{3.1}
T_1\left( {l,m;p} \right): = \sum\limits_{k,n = 1}^\infty  {\frac{{H_k^{\left( l \right)}H_n^{\left( m \right)}}}{{{k^p}n\left( {n + k} \right)}}}\quad{\rm and}\quad T_2\left( {l,m;p} \right): = \sum\limits_{k,n = 1}^\infty  {\frac{{H_n^{\left( l \right)}H_n^{\left( m \right)}}}{{{k^p}n\left( {n + k} \right)}}}\quad \left( {l,m,p \in \mathbb{N}} \right).
\end{align}
In order to prove the main theorems, we need the following lemmas.
\begin{lem}\label{lem3.1}(\cite{X2017})
 Let $m,k$ be positive integers, then
\begin{align}\label{3.2}
\sum\limits_{n = 1}^\infty  {\frac{{{H^{(m)} _n}}}{{n\left( {n + k} \right)}}}
= \frac{1}{k}\left\{ {\zeta( {m + 1})
+ \sum\limits_{j = 1}^{m - 1} {{{\left( { - 1} \right)}^{j - 1}}\zeta( {m + 1 - j}){H^{(j)} _{k - 1}}}  + {{\left( { - 1} \right)}^{m - 1}}\sum\limits_{i = 1}^{k - 1} {\frac{{{H_i}}}{{{i^m}}}} } \right\}.
\end{align}
\end{lem}
\begin{lem}\label{lem3.2}(\cite{SX2017})
For $ {l_1},{l_2},{l_3} \in \N$ and $x,y,z \in \left[ { - 1,1} \right)$, we have the following relation
\begin{align}\label{3.3}
&\sum\limits_{n = 1}^\infty  {\frac{{{\zeta _n}\left( {{l_1};x} \right){\zeta _n}\left( {{l_2};y} \right)}}{{{n^{l_3}}}}{z^n}}  + \sum\limits_{n = 1}^\infty  {\frac{{{\zeta _n}\left( {{l_1};x} \right){\zeta _n}\left( {{l_3};z} \right)}}{{{n^{{l_2}}}}}{y^n}}  + \sum\limits_{n = 1}^\infty  {\frac{{{\zeta _n}\left( {{l_2};y} \right){\zeta _n}\left( {{l_3};z} \right)}}{{{n^{{l_1}}}}}{x^n}}\nonumber \\
& =\sum\limits_{n = 1}^\infty  {\frac{{{\zeta _n}\left( {{l_3};z} \right)}}{{{n^{{l_1} + {l_2}}}}}{{\left( {xy} \right)}^n}}  + \sum\limits_{n = 1}^\infty  {\frac{{{\zeta _n}\left( {{l_1};x} \right)}}{{{n^{{l_3} + {l_2}}}}}{{\left( {yz} \right)}^n}}  + \sum\limits_{n = 1}^\infty  {\frac{{{\zeta _n}\left( {{l_2};y} \right)}}{{{n^{{l_1} + {l_3}}}}}{{\left( {xz} \right)}^n}} \nonumber \\
&\quad  +{\rm Li}{_{l_3}}\left( z \right){\rm Li}{_{{l_1}}}\left( x \right){\rm Li}{_{{l_2}}}\left( y \right) - {\rm Li}{_{{l_1} + {l_2} + {l_3}}}\left( {xyz} \right),
\end{align}
where ${\rm Li}_p$ denotes the polylogarithm function defined by
\[{\rm Li}{_p}\left( x \right) := \sum\limits_{n = 1}^\infty  {\frac{{{x^n}}}{{{n^p}}}}, \Re(p)>1,\ \left| x \right| \le 1 ,\]
with ${\rm Li_1}(x)=-\ln(1-x),\ x\in [-1,1)$.

The partial sum ${\zeta _n}\left( {l;x} \right)$ is defined by \[{\zeta _n}\left( {l;x} \right) := \sum\limits_{k = 1}^n {\frac{{{x^k}}}{k^l}},\quad (l>0,x\in [-1,1]),\]
which is also called the partial sum of polylogarithm function.
\end{lem}
\subsection{Proof of Theorem \ref{thm1}}
From the definition of Tornheim type series $T_1\left( {l,m;p} \right)$, we can rewrite the right hand side of (\ref{3.1}) as
\begin{align}
 T_1\left( {l,m;p} \right) &= \sum\limits_{k = 1}^\infty  {\frac{{H_k^{\left( l \right)}}}{{{k^p}}}\sum\limits_{n = 1}^\infty  {\frac{{H_n^{\left( m \right)}}}{{n\left( {n + k} \right)}}} }\nonumber  \\
  &= \sum\limits_{n = 1}^\infty  {\frac{{H_n^{\left( m \right)}}}{n}} \sum\limits_{k = 1}^\infty  {\frac{{H_k^{\left( l \right)}}}{{{k^p}\left( {n + k} \right)}}}.
\end{align}
Then with the help of the following partial fraction decomposition formula
\[\frac{1}{{{k^p}\left( {n + k} \right)}} = \sum\limits_{i = 1}^{p - 1} {\frac{{{{\left( { - 1} \right)}^{i - 1}}}}{{{n^i}}}\frac{1}{{{k^{p + 1 - i}}}}}  + \frac{{{{\left( { - 1} \right)}^{p - 1}}}}{{{n^{p - 1}}}}\frac{1}{{k\left( {n + k} \right)}},\]
and using Lemma \ref{lem3.1}, we deduce that
\begin{align}\label{3.5}
 &{\left( { - 1} \right)^{m - 1}}\sum\limits_{n = 1}^\infty  {\frac{{H_n^{\left( l \right)}}}{{{n^{p + 1}}}}\left( {\sum\limits_{k = 1}^n {\frac{{{H_k}}}{{{k^m}}}} } \right)}  - {\left( { - 1} \right)^{p + l}}\sum\limits_{n = 1}^\infty  {\frac{{H_n^{\left( m \right)}}}{{{n^{p + 1}}}}\left( {\sum\limits_{k = 1}^n {\frac{{{H_k}}}{{{k^l}}}} } \right)}  \nonumber \\
&= \sum\limits_{i = 1}^{p - 1} {{{\left( { - 1} \right)}^{i - 1}}{S_{m,i + 1}} \cdot {S_{l,p + 1 - i}}}  + {\left( { - 1} \right)^{p - 1}}\zeta \left( {l + 1} \right){S_{m,p + 1}} \nonumber \\
&\quad + {\left( { - 1} \right)^{p - 1}}\sum\limits_{j = 1}^{l - 1} {{{\left( { - 1} \right)}^{j - 1}}\zeta \left( {l + 1 - j} \right)\left\{ {{S_{jm,p + 1}} - {S_{m,p + j + 1}}} \right\}} \nonumber \\
&\quad - {\left( { - 1} \right)^{p + l}}{S_{1m,p + l + 1}} + {\left( { - 1} \right)^{m - 1}}{S_{1l,p + m + 1}} - \zeta \left( {m + 1} \right){S_{l,p + 1}}\nonumber \\
&\quad - \sum\limits_{j = 1}^{m - 1} {{{\left( { - 1} \right)}^{j - 1}}\zeta \left( {m + 1 - j} \right)\left\{ {{S_{jl,p + 1}} - {S_{l,p + j + 1}}} \right\}} .
\end{align}
Note that by a simple calculation, we obtain a simple reflection formula
\begin{align}\label{3.6}
\sum\limits_{k = 1}^n {{A_k}{b_k}}  + \sum\limits_{k = 1}^n {{B_k}{a_k}}  = \sum\limits_{k = 1}^n {{a_k}{b_k}}  + {A_n}{B_n},
\end{align}
where \[{A_n} = \sum\limits_{k = 1}^n {{a_k}} ,{B_n} = \sum\limits_{k = 1}^n {{b_k}}. \]
Hence, from (\ref{3.6}), the following identities are easily derived
\begin{align}\label{3.7}
\sum\limits_{k = 1}^n {\frac{{{H_k}}}{k}}  = \frac{1}{2}\left\{ {H_n^2 + H_n^{\left( 2 \right)}} \right\},
\end{align}
\begin{align}\label{3.8}
\sum\limits_{n = 1}^\infty  {\frac{{{H_n}}}{{{n^p}}}} \left( {\sum\limits_{k = 1}^n {\frac{{{H_k}}}{{{k^m}}}} } \right) + \sum\limits_{n = 1}^\infty  {\frac{{{H_n}}}{{{n^m}}}} \left( {\sum\limits_{k = 1}^n {\frac{{{H_k}}}{{{k^p}}}} } \right) = {S_{1,p}}{S_{1,m}} + {S_{{1^2},p + m}},\;2 \le m,p \in \N.
\end{align}
Thus, letting $l=1$ and replacing $p$ by $p-1$ in (\ref{3.5}), then combining formulas (\ref{3.7}) with (\ref{3.8}), by a direct calculation, we may easily deduce the result (\ref{2.1}). \hfill$\square$

\subsection{Proof of Theorem \ref{thm2}}
Letting $x,y\rightarrow 1$ in (\ref{3.3}), then multiplying it by $(1-z)^{-1}$ and integrating over the interval $(0,1)$, we can find that
\begin{align}\label{3.9}
 &\sum\limits_{n = 1}^\infty  {\frac{{{H_n}H_n^{\left( {{l_1}} \right)}H_n^{\left( {{l_2}} \right)}}}{{{n^{{l_3}}}}} + \sum\limits_{n = 1}^\infty  {\frac{{H_n^{\left( {{l_1}} \right)}\left( {\sum\limits_{k = 1}^n {\frac{{{H_k}}}{{{k^{{l_3}}}}}} } \right)}}{{{n^{{l_2}}}}}}  + \sum\limits_{n = 1}^\infty  {\frac{{H_n^{\left( {{l_2}} \right)}\left( {\sum\limits_{k = 1}^n {\frac{{{H_k}}}{{{k^{{l_3}}}}}} } \right)}}{{{n^{{l_1}}}}}} } \nonumber \\
 & = \sum\limits_{n = 1}^\infty  {\frac{{\left( {\sum\limits_{k = 1}^n {\frac{{{H_k}}}{{{k^{{l_3}}}}}} } \right)}}{{{n^{{l_1} + {l_2}}}}} + \sum\limits_{n = 1}^\infty  {\frac{{{H_n}H_n^{\left( {{l_1}} \right)}}}{{{n^{{l_3} + {l_2}}}}}}  + \sum\limits_{n = 1}^\infty  {\frac{{{H_n}H_n^{\left( {{l_2}} \right)}}}{{{n^{{l_3} + {l_1}}}}}} }  \nonumber \\&\quad + \zeta \left( {{l_1}} \right)\zeta \left( {{l_2}} \right)\left( {\sum\limits_{n = 1}^\infty  {\frac{{{H_n}}}{{{n^{{l_3}}}}}} } \right) - \left( {\sum\limits_{n = 1}^\infty  {\frac{{{H_n}}}{{{n^{{l_3} + {l_1} + {l_2}}}}}} } \right).
\end{align}
From \cite{FS1998,X2017}, we derive the following identity
\begin{align}\label{3.10}
\sum\limits_{n = 1}^\infty  {\frac{1}{{{n^p}}}\left( {\sum\limits_{k = 1}^n {\frac{{{H_k}}}{{{k^m}}}} } \right)}  = \sum\limits_{n = 1}^\infty  {\frac{{{H_n}}}{{{n^{m + p}}}}}  + \zeta \left( p \right)\sum\limits_{n = 1}^\infty  {\frac{{{H_n}}}{{{n^m}}}}  - \sum\limits_{n = 1}^\infty  {\frac{{{H_n}H_n^{\left( p \right)}}}{{{n^m}}}} .
\end{align}
Substituting (\ref{3.10}) into (\ref{3.9}), then
\begin{align}\label{3.11}
 &{S_{1{l_1}{l_2},{l_3}}} + \sum\limits_{n = 1}^\infty  {\frac{{H_n^{\left( {{l_1}} \right)}}}{{{n^{{l_2}}}}}\left( {\sum\limits_{k = 1}^n {\frac{{{H_k}}}{{{k^{{l_3}}}}}} } \right)}  + \sum\limits_{n = 1}^\infty  {\frac{{H_n^{\left( {{l_2}} \right)}}}{{{n^{{l_1}}}}}\left( {\sum\limits_{k = 1}^n {\frac{{{H_k}}}{{{k^{{l_3}}}}}} } \right)} \nonumber \\
  &= {S_{1{l_1},{l_2} + {l_3}}} + {S_{1{l_2},{l_1} + {l_3}}} - {S_{1\left( {{l_1} + {l_2}} \right),{l_3}}} + \zeta \left( {{l_1} + {l_2}} \right){S_{1,{l_3}}} + \zeta \left( {{l_1}} \right)\zeta \left( {{l_2}}. \right){S_{1,{l_3}}}
  \end{align}
Hence, we can easily deduce from (\ref{3.11}) that
\begin{align}\label{3.12}
 &{S_{1{l_1}{l_2},{l_3}}} + {S_{1{l_1}{l_3},{l_2}}} + {S_{1{l_2}{l_3},{l_1}}} \nonumber \\
  &+ \sum\limits_{n = 1}^\infty  {\left\{ {\frac{{H_n^{\left( {{l_1}} \right)}}}{{{n^{{l_2}}}}}\left( {\sum\limits_{k = 1}^n {\frac{{{H_k}}}{{{k^{{l_3}}}}}} } \right) + \frac{{H_n^{\left( {{l_3}} \right)}}}{{{n^{{l_2}}}}}\left( {\sum\limits_{k = 1}^n {\frac{{{H_k}}}{{{k^{{l_1}}}}}} } \right)} \right\}} \nonumber  \\
  &+ \sum\limits_{n = 1}^\infty  {\left\{ {\frac{{H_n^{\left( {{l_2}} \right)}}}{{{n^{{l_1}}}}}\left( {\sum\limits_{k = 1}^n {\frac{{{H_k}}}{{{k^{{l_3}}}}}} } \right) + \frac{{H_n^{\left( {{l_1}} \right)}}}{{{n^{{l_1}}}}}\left( {\sum\limits_{k = 1}^n {\frac{{{H_k}}}{{{k^{{l_2}}}}}} } \right)} \right\}} \nonumber  \\
  &+ \sum\limits_{n = 1}^\infty  {\left\{ {\frac{{H_n^{\left( {{l_1}} \right)}}}{{{n^{{l_3}}}}}\left( {\sum\limits_{k = 1}^n {\frac{{{H_k}}}{{{k^{{l_2}}}}}} } \right) + \frac{{H_n^{\left( {{l_2}} \right)}}}{{{n^{{l_3}}}}}\left( {\sum\limits_{k = 1}^n {\frac{{{H_k}}}{{{k^{{l_1}}}}}} } \right)} \right\}} \nonumber  \\
  &= 2{S_{1{l_1},{l_2} + {l_3}}} + 2{S_{1{l_2},{l_1} + {l_3}}} + 2{S_{1{l_3},{l_1} + {l_2}}} \nonumber \\
 &\quad - {S_{1\left( {{l_1} + {l_2}} \right),{l_3}}} - {S_{1\left( {{l_1} + {l_3}} \right),{l_2}}} - {S_{1\left( {{l_2} + {l_3}} \right),{l_1}}} \nonumber \\
  &\quad+ \left( {\zeta \left( {{l_1} + {l_2}} \right) + \zeta \left( {{l_1}} \right)\zeta \left( {{l_2}} \right)} \right){S_{1,{l_3}}}\nonumber  \\
 &\quad + \left( {\zeta \left( {{l_1} + {l_3}} \right) + \zeta \left( {{l_1}} \right)\zeta \left( {{l_3}} \right)} \right){S_{1,{l_2}}} \nonumber \\
 & \quad+ \left( {\zeta \left( {{l_2} + {l_3}} \right) + \zeta \left( {{l_2}} \right)\zeta \left( {{l_3}} \right)} \right){S_{1,{l_1}}} .
\end{align}
Letting $l=l_1,p=l_2-1,m=l_3$ and $l_1+l_2+l_3=2r+1\ (r\in\N)$ with $l_1,l_2,l_3\geq 2$ in (\ref{3.5}), we obtain
\begin{align}\label{3.13}
 &\sum\limits_{n = 1}^\infty  {\left\{ {\frac{{H_n^{\left( {{l_1}} \right)}}}{{{n^{{l_2}}}}}\left( {\sum\limits_{k = 1}^n {\frac{{{H_k}}}{{{k^{{l_3}}}}}} } \right) + \frac{{H_n^{\left( {{l_3}} \right)}}}{{{n^{{l_2}}}}}\left( {\sum\limits_{k = 1}^n {\frac{{{H_k}}}{{{k^{{l_1}}}}}} } \right)} \right\}} \nonumber \\
  &= {\left( { - 1} \right)^{{l_3} - 1}}\sum\limits_{i = 1}^{{l_2} - 2} {{{\left( { - 1} \right)}^{i - 1}}{S_{{l_3},i + 1}} \cdot {S_{{l_1},{l_2} - i}}}  \nonumber \\
  &\quad+ {\left( { - 1} \right)^{{l_2} + {l_3} - 1}}\sum\limits_{j = 1}^{{l_1} - 1} {{{\left( { - 1} \right)}^{j - 1}}\zeta \left( {{l_1} + 1 - j} \right)\left( {{S_{{l_3}j,{l_2}}} - {S_{{l_3},{l_2} + j}}} \right)}  \nonumber \\
  &\quad- {\left( { - 1} \right)^{{l_3} - 1}}\sum\limits_{j = 1}^{{l_3} - 1} {{{\left( { - 1} \right)}^{j - 1}}\zeta \left( {{l_3} + 1 - j} \right)\left( {{S_{{l_1}j,{l_2}}} - {S_{{l_1},{l_2} + j}}} \right)}  \nonumber \\
  &\quad+ {S_{1{l_3},{l_1} + {l_2}}} + {S_{1{l_1},{l_2} + {l_3}}} - {\left( { - 1} \right)^{{l_3} - 1}}{S_{{l_1},{l_2}}} + {\left( { - 1} \right)^{{l_2} + {l_3} - 1}}\zeta \left( {{l_1} + 1} \right){S_{{l_3},{l_2}}}.
\end{align}
Thus, combining formulas (\ref{3.12}) and (\ref{3.13}) and noting that the symmetry relations
\[{S_{p,q}} + {S_{q,p}} = \zeta \left( {p + q} \right) + \zeta \left( p \right)\zeta \left( q \right),\quad (p,q>1),\]
we deduce the desired result.  This completes the proof of Theorem \ref{thm2}.\hfill$\square$
\subsection{Proof of Theorem \ref{thm3}}
In this subsection, we consider the Tornheim type series $T_2\left( {m_1,m_2;p} \right)$ to prove the Theorem \ref{thm3}. First, from \cite{XMZ2016}, we have
\begin{align}\label{3.14}
\sum\limits_{n = 1}^\infty  {\frac{{H_n^{\left( {{m_1}} \right)}H_n^{\left( {{m_2}} \right)}}}{{n\left( {n + k} \right)}}}  =& \frac{1}{k}\sum\limits_{n = 1}^\infty  {\left\{ {\frac{{H_n^{\left( {{m_1}} \right)}}}{{{n^{{m_2} + 1}}}} + \frac{{H_n^{\left( {{m_2}} \right)}}}{{{n^{{m_1} + 1}}}}} \right\}}  - \sum\limits_{n = 1}^\infty  {\frac{{H_n^{\left( {{m_1} + {m_2}} \right)}}}{{n\left( {n + k} \right)}}}\nonumber \\
& + \frac{1}{k}\sum\limits_{j = 1}^{k - 1} {\sum\limits_{n = 1}^\infty  {\left\{ {\frac{{H_n^{\left( {{m_1}} \right)}}}{{{n^{{m_2}}}\left( {n + j} \right)}} + \frac{{H_n^{\left( {{m_2}} \right)}}}{{{n^{{m_1}}}\left( {n + j} \right)}}} \right\}} },\quad m_1,m_2,k\in\N.
\end{align}
Then, applying the formula (\ref{3.2}) into (\ref{3.14}), by a simple calculation we obtain the following result
\begin{align}\label{3.15}
\sum\limits_{n = 1}^\infty  {\frac{{H_n^{\left( {{m_1}} \right)}H_n^{\left( {{m_2}} \right)}}}{{n\left( {n + k} \right)}}}  = \frac{1}{k}\left\{ \begin{array}{l}
 {S_{{m_1},{m_2} + 1}} + {S_{{m_2},{m_1} + 1}} - \zeta \left( {{m_1} + {m_2} + 1} \right) \\
  + \sum\limits_{i = 1}^{{m_2} - 1} {{{\left( { - 1} \right)}^{i - 1}}{S_{{m_1},{m_2} + 1 - i}}{\zeta _{k - 1}}\left( i \right)}  \\
  + \sum\limits_{i = 1}^{{m_1} - 1} {{{\left( { - 1} \right)}^{i - 1}}{S_{{m_2},{m_1} + 1 - i}}{\zeta _{k - 1}}\left( i \right)}  \\
  + {\left( { - 1} \right)^{{m_2} - 1}}\zeta \left( {{m_1} + 1} \right){\zeta _{k - 1}}\left( {{m_2}} \right) \\
  + {\left( { - 1} \right)^{{m_1} - 1}}\zeta \left( {{m_2} + 1} \right){\zeta _{k - 1}}\left( {{m_1}} \right) \\
  + {\left( { - 1} \right)^{{m_2} - 1}}\sum\limits_{i = 1}^{{m_1} - 1} {{{\left( { - 1} \right)}^{i - 1}}\zeta \left( {{m_1} + 1 - i} \right){\zeta _{k - 1}}\left( {{m_2},i} \right)}  \\
  + {\left( { - 1} \right)^{{m_1} - 1}}\sum\limits_{i = 1}^{{m_2} - 1} {{{\left( { - 1} \right)}^{i - 1}}\zeta \left( {{m_2} + 1 - i} \right){\zeta _{k - 1}}\left( {{m_1},i} \right)}  \\
  + {\left( { - 1} \right)^{{m_1} + {m_2}}}\left( {{\zeta _{k - 1}}\left( {{m_2},{m_1},1} \right) + {\zeta _{k - 1}}\left( {{m_2},{m_1} + 1} \right)} \right) \\
  + {\left( { - 1} \right)^{{m_1} + {m_2}}}\left( {{\zeta _{k - 1}}\left( {{m_1},{m_2},1} \right) + {\zeta _{k - 1}}\left( {{m_1},{m_2} + 1} \right)} \right) \\
  + {\left( { - 1} \right)^{{m_1} + {m_2}}}\left( {{\zeta _{k - 1}}\left( {{m_1} + {m_2},1} \right) + {\zeta _{k - 1}}\left( {{m_1} + {m_2} + 1} \right)} \right) \\
  - \sum\limits_{i = 1}^{{m_1} + {m_2} - 1} {{{\left( { - 1} \right)}^{i - 1}}\zeta \left( {{m_1} + {m_2} + 1 - i} \right){\zeta _{k - 1}}\left( i \right)} , \\
 \end{array} \right\},
 \end{align}
where ${\zeta _n}\left( {{s_1},{s_2}, \cdots ,{s_k}} \right)$ denotes the multiple harmonic sum (also called the partial sum of multiple zeta value) defined by \cite{Xu2017,KP2013}
\[{\zeta _n}\left( {{s_1},{s_2}, \cdots ,{s_k}} \right): = \sum\limits_{n \ge {n_1} > {n_2} >  \cdots  > {n_k} \ge 1} {\frac{1}{{n_1^{{s_1}}n_2^{{s_2}} \cdots n_k^{{s_k}}}}} ,\]
when $n<k$, then ${\zeta _n}\left( {{s_1},{s_2}, \cdots ,{s_k}} \right)=0$, and ${\zeta _n}\left(\emptyset \right)=1$. The integers $k$ and $w:=s_1+\ldots+s_k$ are called the depth and the weight of a multiple harmonic sum. It is obvious that \[{\zeta _n}\left( s \right) = H_n^{\left( s \right)},\quad s\in\N.\]

By the definition of $T_2\left( {m_1,m_2;p} \right)$ we can find that
\begin{align*}
{T_2}\left( {{m_1},{m_2};p} \right) &= \sum\limits_{k = 1}^\infty  {\frac{1}{{{k^p}}}\sum\limits_{n = 1}^\infty  {\frac{{H_n^{\left( {{m_1}} \right)}H_n^{\left( {{m_2}} \right)}}}{{n\left( {n + k} \right)}}} } \\&= \sum\limits_{n = 1}^\infty  {\frac{{H_n^{\left( {{m_1}} \right)}H_n^{\left( {{m_2}} \right)}}}{n}} \sum\limits_{k = 1}^\infty  {\frac{1}{{{k^p}\left( {n + k} \right)}}}.
\end{align*}

Hence, multiplying (\ref{3.15}) by $\frac 1{k^p}$ and taking the summation from $k=0$ to $\infty$, by a similar argument as in the proofs of Theorem \ref{thm1} and \ref{thm2}, we may easily deduce the dedired result. \hfill$\square$

Taking $m_1=m_2=1,m_3=2$ in (\ref{2.3}) yields
\[{S_{{1^3},2}} = 4\zeta \left( 2 \right)\zeta \left( 3 \right) + 2\zeta \left( {2,{{\left\{ 1 \right\}}_3}} \right) + 2\zeta \left( {2,1,2} \right) + \zeta \left( {2,2,1} \right) + \zeta \left( {2,3} \right).\]

\section{Some Examples and Corollaries}
From Theorem \ref{thm1}, \ref{thm2} and Corollary \ref{cor3}, \ref{cor4} in section 2, we can get the following some connections between Euler sums and zeta values.
\begin{align}
&{S_{{1^2}2,2}} = \frac{{41}}{{12}}\zeta \left( 6 \right) + 2{\zeta ^2}\left( 3 \right),\\
&{S_{{1^2}3,3}} = \frac{9}{2}\zeta \left( 3 \right)\zeta \left( 5 \right) + \frac{3}{2}\zeta \left( 2 \right){\zeta ^2}\left( 3 \right) - \frac{{443}}{{288}}\zeta \left( 8 \right) - \frac{{23}}{4}{S_{2,6}},\\
&{S_{{{13}^2},3}} = \frac{{883}}{{20}}\zeta \left( {10} \right) - 26{\zeta ^2}\left( 5 \right) - \frac{{31}}{4}\zeta \left( 3 \right)\zeta \left( 7 \right) - 8\zeta \left( 2 \right)\zeta \left( 3 \right)\zeta \left( 5 \right)\nonumber\\
&\quad\quad\quad\quad + \frac{3}{4}{\zeta ^2}\left( 3 \right)\zeta \left( 4 \right) + 9\zeta \left( 2 \right){S_{2,6}} - \frac{{21}}{4}{S_{2,8}},\\
&{S_{{1^2}4,4}} =  \frac{{7749}}{{160}}\zeta \left( {10} \right) - 16{\zeta ^2}\left( 5 \right) - \frac{{125}}{8}\zeta \left( 3 \right)\zeta \left( 7 \right) - 14\zeta \left( 2 \right)\zeta \left( 3 \right)\zeta \left( 5 \right)\nonumber\\
&\quad\quad\quad\quad + \frac{3}{2}{\zeta ^2}\left( 3 \right)\zeta \left( 4 \right) + 11{S_{2,8}} + \frac{5}{2}\zeta \left( 2 \right){S_{2,6}},\\
&{S_{{{12}^2},3}} + 2{S_{123,2}} =  - \frac{{2225}}{{96}}\zeta \left( 8 \right) + \frac{{73}}{2}\zeta \left( 3 \right)\zeta \left( 5 \right) - \frac{5}{2}\zeta \left( 2 \right){\zeta ^2}\left( 3 \right) - \frac{{31}}{4}{S_{2,6}},\\
&{S_{{1^2}2,5}} + {S_{{1^2}5,2}} = \frac{{1235}}{{72}}\zeta \left( 9 \right) - 2{\zeta ^3}\left( 3 \right) - \frac{{101}}{8}\zeta (2)\zeta (7) + \frac{{185}}{{24}}\zeta \left( 3 \right)\zeta \left( 6 \right) + \frac{{13}}{4}\zeta (4)\zeta (5),\\
&{S_{{1^2}3,4}} + {S_{{1^2}4,3}} = \frac{{761}}{{36}}\zeta \left( 9 \right) + {\zeta ^3}\left( 3 \right) - \frac{3}{8}\zeta \left( 2 \right)\zeta \left( 7 \right) - 11\zeta \left( 3 \right)\zeta \left( 6 \right) - \frac{{21}}{4}\zeta \left( 4 \right)\zeta \left( 5 \right),\\
&{S_{{{12}^2},5}} + 2{S_{125,2}} =  - \frac{{2403}}{{160}}\zeta \left( {10} \right) + \frac{{69}}{4}{\zeta ^2}\left( 5 \right) + \frac{{491}}{8}\zeta \left( 3 \right)\zeta \left( 7 \right) - 34\zeta \left( 2 \right)\zeta \left( 3 \right)\zeta \left( 5 \right)\nonumber\\
&\quad\quad\quad\quad\quad\quad\quad\quad - \frac{{27}}{4}{\zeta ^2}\left( 3 \right)\zeta \left( 4 \right) + \frac{{27}}{2}\zeta \left( 2 \right){S_{2,6}} - 13{S_{2,8}},\\
&{S_{123,4}} + {S_{124,3}} + {S_{134,2}} =  - \frac{{12889}}{{320}}\zeta \left( {10} \right) + \frac{{85}}{4}{\zeta ^2}\left( 5 \right) + \frac{{151}}{{16}}\zeta \left( 3 \right)\zeta \left( 7 \right) + \frac{{25}}{2}\zeta \left( 2 \right)\zeta \left( 3 \right)\zeta \left( 5 \right)\nonumber\\
&\quad\quad\quad\quad\quad\quad\quad\quad\quad\quad\quad - \frac{9}{8}{\zeta ^2}\left( 3 \right)\zeta \left( 4 \right) - \frac{{39}}{4}\zeta \left( 2 \right){S_{2,6}} + 4{S_{2,8}},\\
&{S_{{1^2}2,6}} + {S_{{1^2}6,2}} + {S_{{2^2},6}} + {S_{26,2}} = \frac{{5837}}{{160}}\zeta \left( {10} \right) - 20{\zeta ^2}\left( 5 \right) - \frac{{189}}{8}\zeta \left( 3 \right)\zeta \left( 7 \right) + 7\zeta \left( 2 \right)\zeta \left( 3 \right)\zeta \left( 5 \right)\nonumber\\
& \quad\quad\quad\quad\quad\quad\quad\quad\quad\quad\quad\quad\quad\quad+ \frac{1}{4}{\zeta ^2}\left( 3 \right)\zeta \left( 4 \right) - \frac{3}{2}\zeta \left( 2 \right){S_{2,6}} + 10{S_{2,8}},\\
&{S_{{1^2}3,5}} + {S_{{1^2}5,3}} + {S_{23,5}} + {S_{25,3}} =  - \frac{{4559}}{{60}}\zeta \left( {10} \right) + \frac{{103}}{2}{\zeta ^2}\left( 5 \right) + \frac{{337}}{4}\zeta \left( 3 \right)\zeta \left( 7 \right) - 29\zeta \left( 2 \right)\zeta \left( 3 \right)\zeta \left( 5 \right)\nonumber\\
&  \quad\quad\quad\quad\quad\quad\quad\quad\quad\quad\quad\quad\quad\quad- \frac{3}{2}{\zeta ^2}\left( 3 \right)\zeta \left( 4 \right) + 9\zeta \left( 2 \right){S_{2,6}} - \frac{{59}}{2}{S_{2,8}}.
\end{align}
It is known that the quadratic Euler sums of even weight like ${S_{{p_1}{p_2},q}}: = \sum\limits_{n = 1}^\infty  {H_n^{\left( {{p_1}} \right)}H_n^{\left( {{p_2}} \right)}{n^{ - q}}} $ can be expressed by linear sums and zeta values at integer arguments if $p_1+p_2+q$ is even and $p_1,p_2,q>1$, see \cite{FS1998}. We can deduce the following results
\begin{align}
S_{2^2,6}  =& \frac{{2697}}
{{40}}\zeta( {10}) - 41{\zeta ^2}( 5) - 63\zeta( 3)\zeta( 7) + 16\zeta( 2)\zeta( 3)\zeta( 5)\nonumber\\
& + 4{\zeta ^2}( 3)\zeta( 4) + \frac{{23}}
{2}{S_{2,8}}{\text{ + }}2\zeta( 2){S_{2,6}},\\
{S_{26,2}} =&  - \frac{{2997}}{{80}}\zeta \left( {10} \right) + 23{\zeta ^2}\left( 5 \right) + 35\zeta \left( 3 \right)\zeta \left( 7 \right) - 8\zeta \left( 2 \right)\zeta \left( 3 \right)\zeta \left( 5 \right)\nonumber\\
& - 2{\zeta ^2}\left( 3 \right)\zeta \left( 4 \right) - \frac{{13}}{2}{S_{2,8}} - \zeta \left( 2 \right){S_{2,6}},\\
{S_{23,5}}  = & - \frac{{2227}}{{32}}\zeta \left( {10} \right) + \frac{{89}}{2}{\zeta ^2}\left( 5 \right) + 56\zeta \left( 3 \right)\zeta \left( 7 \right) - 15\zeta \left( 2 \right)\zeta \left( 3 \right)\zeta \left( 5 \right)\nonumber\\
& - 10{S_{2,8}} - \frac{5}{2}\zeta \left( 2 \right){S_{2,6}},\\
{S_{25,3}}  =& \frac{{3223}}{{160}}\zeta \left( {10} \right) - \frac{{17}}{4}{\zeta ^2}\left( 5 \right) + \frac{{63}}{2}\zeta \left( 3 \right)\zeta \left( 7 \right) - 32\zeta \left( 2 \right)\zeta \left( 3 \right)\zeta \left( 5 \right)\nonumber\\
& - 2{\zeta ^2}\left( 3 \right)\zeta \left( 4 \right) - \frac{{19}}{4}{S_{2,8}}{\rm{ + }}\frac{{25}}{2}\zeta \left( 2 \right){S_{2,6}}.
\end{align}
Hence, we obtain the following closed form of two combined sums involving two cubic Euler sums
\begin{align}
{S_{{1^2}2,6}} + {S_{{1^2}6,2}} =& \frac{{1043}}
{{160}}\zeta \left( {10} \right) - 2{\zeta ^2}\left( 5 \right) + \frac{{35}}
{8}\zeta \left( 3 \right)\zeta \left( 7 \right) - \zeta \left( 2 \right)\zeta \left( 3 \right)\zeta \left( 5 \right)\nonumber\\
&- \frac{7}{4}{\zeta ^2}\left( 3 \right)\zeta \left( 4 \right) - \frac{5}
{2}\zeta \left( 2 \right){S_{2,6}} + 5{S_{2,8}},\\
{S_{{1^2}3,5}} + {S_{{1^2}5,3}} =&  - \frac{{398}}
{{15}}\zeta \left( {10} \right) + \frac{{45}}
{4}{\zeta ^2}\left( 5 \right) - \frac{{13}}
{4}\zeta \left( 3 \right)\zeta \left( 7 \right) + 18\zeta \left( 2 \right)\zeta \left( 3 \right)\zeta \left( 5 \right) \nonumber\\&+ \frac{1}
{2}{\zeta ^2}\left( 3 \right)\zeta \left( 4 \right) - \zeta \left( 2 \right){S_{2,6}} - \frac{{59}}
{4}{S_{2,8}}.
\end{align}
In \cite{Xu2017}, we prove that for integers $m,p> 0$ and $r>1$, then
\begin{align}\label{4.12}
\sum\limits_{n = 1}^\infty  {\frac{{s\left( {n,m} \right){H^{(r)}_n}}}{{n!{n^p}}}}  = \zeta \left( r \right)\zeta \left( {m + 1,{{\left\{ 1 \right\}}_{p - 1}}} \right) - \zeta \left( {m + 1,{{\left\{ 1 \right\}}_{p - 1}},2,{{\left\{ 1 \right\}}_{r - 2}}} \right).
\end{align}
Here  ${s\left( {n,k} \right)}$ denotes the (unsigned) Stirling number of the first kind \cite{C1974}, and we have
\begin{align*}
& s\left( {n,1} \right) = \left( {n - 1} \right)!,\\&s\left( {n,2} \right) = \left( {n - 1} \right)!{H_{n - 1}},\\&s\left( {n,3} \right) = \frac{{\left( {n - 1} \right)!}}{2}\left[ {H_{n - 1}^2 - {H^{(2)} _{n - 1}}} \right],\\
&s\left( {n,4} \right) = \frac{{\left( {n - 1} \right)!}}{6}\left[ {H_{n - 1}^3 - 3{H_{n - 1}}{H^{(2)} _{n - 1}} + 2{H^{(3)}_{n - 1}}} \right], \\
&s\left( {n,5} \right) = \frac{{\left( {n - 1} \right)!}}{{24}}\left[ {H_{n - 1}^4 - 6{H^{(4)}_{n - 1}} - 6H_{n - 1}^2{H^{(2)}_{n - 1}} + 3(H^{(2)}_{n-1})^2+ 8H_{n - 1}^{}{H^{(3)}_{n - 1}}} \right].
\end{align*}
Hence, taking $m=3$ in (\ref{4.12}) and using the duality relation (\ref{1.6}), we obtain
\begin{align}\label{4.13}
{S_{{1^2}r,p + 1}} = {S_{2r,p + 1}} + 2{S_{1r,p + 2}} - 2{S_{r,p + 3}} + 2\zeta \left( r \right)\zeta \left( {4,{{\left\{ 1 \right\}}_{p - 1}}} \right) - 2\zeta \left( {r,p + 1,1,1} \right).
\end{align}
Since the multiple zeta value $\zeta(m+1,\{1\}_{k-1} )$ can be
represented as a polynomial of zeta values with rational coefficients \cite{Xu2017}. For example:
\[\begin{array}{l}
 \zeta \left( {2,{{\left\{ 1 \right\}}_m}} \right) = \zeta \left( {m + 2} \right), \\
 \zeta \left( {3,{{\left\{ 1 \right\}}_m}} \right) = \frac{{m + 2}}{2}\zeta \left( {m + 3} \right) - \frac{1}{2}\sum\limits_{k = 1}^m {\zeta \left( {k + 1} \right)\zeta \left( {m + 2 - k} \right)} . \\
 \end{array}\]
Thus, from Theorem \ref{thm1}, \ref{thm2} and \ref{thm3} with the help of the formulas (\ref{1.5}) and (\ref{4.13}), we obtain the following description of multiple zeta values of depth four.
\begin{cor}\label{cor4.1} For positive integers $m_1,m_2,m_3>1$, then the combined sums
\[\zeta \left( {{m_1},{m_2},1,1} \right) + \zeta \left( {{m_2},{m_1},1,1} \right),\quad \zeta \left( {{m_1},1,{m_2},1} \right) + \zeta \left( {{m_2},1,{m_1},1} \right)\]
and
\begin{align*}
&\zeta \left( {{m_1},{m_2},{m_3},1} \right) + \zeta \left( {{m_1},{m_3},{m_2},1} \right) + \zeta \left( {{m_2},{m_3},{m_1},1} \right)\\& + \zeta \left( {{m_2},{m_1},{m_3},1} \right) + \zeta \left( {{m_3},{m_2},{m_1},1} \right) + \zeta \left( {{m_3},{m_1},{m_2},1} \right)
\end{align*}
can be expressed as a rational linear combination of products of single, double or triple zeta values.
\end{cor}

Hence, from Corollary \ref{cor4.1}, we know that for any positive integer $m>1$, then the three quadruple zeta values \[\zeta \left( {m,m,1,1} \right),\zeta \left( {m,1,m,1} \right)\ {\rm and}\ \zeta \left( {m,m,m,1} \right)\]
can be expressed in terms of other multiple zeta values of depth $\leq 3$. Note that the conclusion was also proved in \cite{EW2012} by another method.

On the other hand, we note that H.N. Minh and M. Petitot \cite{M2000} gave all closed form of multiple zeta values of weight $\leq 9$. By using their results and combining (\ref{4.13}), we give the following explicit formulas of cubic Euler sums
\[\begin{array}{l}
 {S_{{1^2}3,4}} = \sum\limits_{n = 1}^\infty  {\frac{{H_n^2H_n^{\left( 3 \right)}}}{{{n^4}}}}  = \frac{{3895}}{{72}}\zeta (9) - \frac{5}{8}\zeta (2)\zeta (7) - \frac{{227}}{{24}}\zeta (3)\zeta (6) - \frac{{75}}{2}\zeta (4)\zeta (5) + {\zeta ^3}(3), \\
 {S_{{1^2}4,3}} = \sum\limits_{n = 1}^\infty  {\frac{{H_n^2H_n^{\left( 4 \right)}}}{{{n^3}}}}  =  - \frac{{791}}{{24}}\zeta (9) + \frac{1}{4}\zeta (2)\zeta (7) - \frac{{37}}{{24}}\zeta (3)\zeta (6) + \frac{{129}}{4}\zeta (4)\zeta (5), \\
 {S_{{1^2}5,2}} = \sum\limits_{n = 1}^\infty  {\frac{{H_n^2H_n^{\left( 5 \right)}}}{{{n^2}}}}  = \frac{{679}}{{18}}\zeta (9) - \frac{{61}}{8}\zeta (2)\zeta (7) - \frac{{55}}{{12}}\zeta (3)\zeta (6) - \frac{{59}}{4}\zeta (4)\zeta (5) + {\zeta ^3}(3). \\
 \end{array}\]

\section{Alternating Cubic Euler Sums}
Let
\[\bar H_n^{\left( p \right)}: = \sum\limits_{j = 1}^n {\frac{{{{\left( { - 1} \right)}^{j - 1}}}}{{{j^p}}}} ,\quad {{\bar H}_n}: = \bar H_n^{\left( 1 \right)}\]
denote the alternating harmonic numbers, $p\in \N$. In this section, we will prove that the alternating cubic Euler sums
\[{S_{{{\bar 1}^3},2r + 1}}: = \sum\limits_{n = 1}^\infty  {\frac{{\bar H_n^3}}{{{n^{2r + 1}}}}} ,\;\left( {r \in \N} \right)\]
are reducible to alternating quadratic and linear Euler sums. The generalized alternating Euler sums are defined by
\[{S_{{P_l}{{\bar R}_k},q}}: = \sum\limits_{n = 1}^\infty  {\frac{{H_n^{\left( {{p_1}} \right)} \cdots H_n^{\left( {{p_l}} \right)}\bar H_n^{\left( {{r_1}} \right)} \cdots \bar H_n^{\left( {{r_k}} \right)}}}{{{n^q}}}} ,\quad {S_{{P_l}{{\bar R}_k},\bar q}}: = \sum\limits_{n = 1}^\infty  {\frac{{H_n^{\left( {{p_1}} \right)} \cdots H_n^{\left( {{p_l}} \right)}\bar H_n^{\left( {{r_1}} \right)} \cdots \bar H_n^{\left( {{r_k}} \right)}}}{{{n^q}}}{{\left( { - 1} \right)}^{n - 1}}}, \]
where the ${P_l}: = \left( {{p_1}, \ldots ,{p_l}} \right)\in (\N_0)^l,{{R}_k}: = \left( {{r_1}, \ldots ,{r_k}} \right)\in (\N_0)^k$, ${{\bar R}_k}: = \left( {{\bar r_1}, \ldots ,{\bar r_k}} \right)$ and the quantity $w: = {p_1} +  \cdots  + {p_l} + {r_1} +  \cdots  + {r_k} + q$ being called the weight and the quantity $l+k$ being the degree, $\N_0:=\N\cup \{0\}$. Since repeated summands in partitions are indicated by powers, we denote, for example, the sums
\[{S_{{1^3}2{{\bar 1}^2}\bar 3,5}} = \sum\limits_{n = 1}^\infty  {\frac{{H_n^3H_n^{\left( 2 \right)}\bar H_n^2\bar H_n^{\left( 3 \right)}}}{{{n^5}}}},\quad {S_{{1^2}3{{\bar 1}^3}\bar 3,\bar 4}} = \sum\limits_{n = 1}^\infty  {\frac{{H_n^2H_n^{\left( 3 \right)}\bar H_n^3\bar H_n^{\left( 3 \right)}}}{{{n^4}}}{{\left( { - 1} \right)}^{n - 1}}} .\]
In \cite{FS1998}, P. Flajolet and B. Salvy gave explicit reductions to zeta values for all linear sums
\[{S_{p,q}} = \sum\limits_{n = 1}^\infty  {\frac{{H_n^{\left( p \right)}}}{{{n^q}}}} ,\:{S_{\bar p,q}} = \sum\limits_{n = 1}^\infty  {\frac{{\bar H_n^{\left( p \right)}}}{{{n^q}}}} ,\:{S_{p,\bar q}} = \sum\limits_{n = 1}^\infty  {\frac{{H_n^{\left( p \right)}}}{{{n^q}}}} {\left( { - 1} \right)^{n - 1}},\:{S_{\bar p,\bar q}} = \sum\limits_{n = 1}^\infty  {\frac{{\bar H_n^{\left( p \right)}}}{{{n^q}}}{{\left( { - 1} \right)}^{n - 1}}} \]
when $p+q$ is an odd weight. The evaluation of linear sums in terms of values of Riemann zeta function and polylogarithm function at positive integers is known when $(p,q) = (1,3), (2,2)$, or $p + q$ is odd (see \cite{BBG1995,FS1998,Xu2016,Xyz2016}). In \cite{Xyz2016}, the author jointly with Y. Yang and J. Zhang proved that the alternating quadratic Euler sums
\begin{align*}
 &{S_{{{\bar 1}^2},{\overline {2m}}}} = \sum\limits_{n = 1}^\infty  {\frac{{\bar H_n^2}}{{{n^{2m}}}}} {\left( { - 1} \right)^{n - 1}},{S_{1\bar 1,{\overline {2m}}}} = \sum\limits_{n = 1}^\infty  {\frac{{{H_n}{{\bar H}_n}}}{{{n^{2m}}}}} {\left( { - 1} \right)^{n - 1}}, \\
& {S_{{1^2},\overline {2m}}} = \sum\limits_{n = 1}^\infty  {\frac{{H_n^2}}{{{n^{2m}}}}} {\left( { - 1} \right)^{n - 1}},{S_{1\bar 1,2m}} = \sum\limits_{n = 1}^\infty  {\frac{{{H_n}{{\bar H}_n}}}{{{n^{2m}}}}} ,{S_{{{\bar 1}^2},2m}} = \sum\limits_{n = 1}^\infty  {\frac{{\bar H_n^2}}{{{n^{2m}}}}} .
\end{align*}
are reducible to polynomials in zeta values and to linear sums, and are able to give explicit values for certain of these sums in terms of the Riemann zeta function and the polylogarithm function, here $m$ is a positive integer. In \cite{BBG1994}, D.H. Bailey, J.M. Borwein and R. Girgensohn considered sums
\[{S_{{{\bar 1}^p},q}} = \sum\limits_{n = 1}^\infty  {\frac{{\bar H_n^p}}{{{n^q}}}} ,{S_{{{\bar 1}^p},\bar q}} = \sum\limits_{n = 1}^\infty  {\frac{{\bar H_n^p}}{{{n^q}}}{{\left( { - 1} \right)}^{n - 1}}} ,{S_{{1^p},\bar q}} = \sum\limits_{n = 1}^\infty  {\frac{{H_n^p}}{{{n^q}}}{{\left( { - 1} \right)}^{n - 1}}} \:\]
and obtained a number of experimental identities using the experimental method, where $p$ and $q$ are positive integers with $p+q\leq 5$.

Now we state our main results in this section. First, we need the following lemma.
\begin{lem} For positive integers $m$ and $k$, then
\begin{align}\label{5.1}
\sum\limits_{n = 1}^\infty  {\frac{{\bar H_n^{\left( m \right)}}}{{n\left( {n + k} \right)}}}  = \frac{1}{k}\left\{ {\begin{array}{*{20}{c}}
   {\bar \zeta \left( {m + 1} \right) + \sum\limits_{j = 1}^{m - 1} {{{\left( { - 1} \right)}^{j - 1}}\bar \zeta \left( {m + 1 - j} \right)H_{k - 1}^{\left( j \right)}} }  \\
   { + {{\left( { - 1} \right)}^{m - 1}}\ln 2\left( {H_{k - 1}^{\left( m \right)} + \bar H_{k - 1}^{\left( m \right)}} \right) + {{\left( { - 1} \right)}^m}\sum\limits_{i = 1}^{k - 1} {\frac{{{{\left( { - 1} \right)}^{i - 1}}}}{{{i^m}}}{{\bar H}_i}} }  \\
\end{array}} \right\},
\end{align}
where $\bar \zeta \left( s \right)$ denotes the alternating Riemann zeta function defined by
$$\bar \zeta \left( s \right) := \sum\limits_{n = 1}^\infty  {\frac{{{{\left( { - 1} \right)}^{n - 1}}}}{{{n^s}}}} ,\;{\mathop{\Re}\nolimits} \left( s \right) \ge 1.$$
\end{lem}
\pf By the definition of polylogarithm function and using Cauchy product formula, we can find that
\begin{align}\label{5.2}
- \frac{{{\rm Li}{_m}\left( { - x} \right)}}{{1 - x}}=\sum\limits_{n = 1}^\infty  {{\bar H}^{(m)}_n{x^n}},\ x\in (-1,1).
\end{align}
Multiplying (\ref{5.2}) by $x^{r-1}-x^{k-1}$ and integrating over (0,1), we obtain
\begin{align}\label{5.3}
\int\limits_0^1 {\left( {{x^{k - 1}} - {x^{r - 1}}} \right)\frac{{{\rm Li}{_m}\left( { - x} \right)}}{{1 - x}}} dx = \left( {k - r} \right)\sum\limits_{n = 1}^\infty  {\frac{{{\bar H}^{(m)}_n}}{{\left( {n + r} \right)\left( {n + k} \right)}}} \:\:\:\left( {0 \le r < k,\:r,k \in \N} \right).
\end{align}
We now evaluate the integral on the left side of (\ref{5.3}). Noting that the integral of (\ref{5.3}) can be written as
\begin{align}\label{5.4}
\int\limits_0^1 {\left( {{x^{k - 1}} - {x^{r - 1}}} \right)\frac{{{\rm Li}{_m}\left( { - x} \right)}}{{1 - x}}} dx = \sum\limits_{i = 1}^{k - r} {{{\left( { - 1} \right)}^{r + i}}\int\limits_0^{ - 1} {{x^{r + i - 2}}{\rm Li}{_m}\left( x \right)} } dx.
\end{align}
Then using integration by parts we have
\begin{align}\label{5.5}
\int\limits_0^x {{t^{n - 1}}{\rm Li}{_q}\left( t \right)dt}  = \sum\limits_{i = 1}^{q - 1} {{{\left( { - 1} \right)}^{i - 1}}\frac{{{x^n}}}{{{n^i}}}{\rm Li}{_{q + 1 - i}}\left( x \right)}  + \frac{{{{\left( { - 1} \right)}^q}}}{{{n^q}}}\ln \left( {1 - x} \right)\left( {{x^n} - 1} \right) - \frac{{{{\left( { - 1} \right)}^q}}}{{{n^q}}}\left( {\sum\limits_{k = 1}^n {\frac{{{x^k}}}{k}} } \right).
\end{align}
Thus, taking $r=0$ in (\ref{5.3}) and (\ref{5.4}), then substituting (\ref{5.5}) into (\ref{5.3}) and (\ref{5.4}), we deduce the formula (\ref{5.1}).\hfill$\square$
\begin{thm} For positive integers $l,m$ and $p$, we have
\begin{align}\label{5.6}
& {\left( { - 1} \right)^m}\sum\limits_{n = 1}^\infty  {\frac{{\bar H_n^{\left( l \right)}}}{{{n^{p + 1}}}}\left( {\sum\limits_{i = 1}^n {\frac{{{{\bar H}_i}}}{{{i^m}}}{{\left( { - 1} \right)}^{i - 1}}} } \right)}  + {\left( { - 1} \right)^{p + l}}\sum\limits_{n = 1}^\infty  {\frac{{\bar H_n^{\left( m \right)}}}{{{n^{p + 1}}}}\left( {\sum\limits_{i = 1}^n {\frac{{{{\bar H}_i}}}{{{i^l}}}{{\left( { - 1} \right)}^{i - 1}}} } \right)} \nonumber \\
  &= \sum\limits_{i = 1}^{p - 1} {{{\left( { - 1} \right)}^{i - 1}}{S_{\bar m,i + 1}} \cdot {S_{\bar l,p + 1 - i}}}  + {\left( { - 1} \right)^{p - 1}}\bar \zeta \left( {l + 1} \right){S_{\bar m,p + 1}} - \bar \zeta \left( {m + 1} \right){S_{\bar l,p + 1}}\nonumber \\
 &\quad + {\left( { - 1} \right)^{p - 1}}\sum\limits_{j = 1}^{l - 1} {{{\left( { - 1} \right)}^{j - 1}}\bar \zeta \left( {l + 1 - j} \right)\left( {{S_{\bar mj,p + 1}} - {S_{\bar m,p + j + 1}}} \right)} \nonumber \\
 &\quad - \sum\limits_{j = 1}^{m - 1} {{{\left( { - 1} \right)}^{j - 1}}\bar \zeta \left( {m + 1 - j} \right)\left( {{S_{\bar lj,p + 1}} - {S_{\bar l,p + j + 1}}} \right)} \nonumber \\
  &\quad+ {\left( { - 1} \right)^{p + l}}\ln 2\left( {{S_{\bar ml,p + 1}} + {S_{\bar m\bar l,p + 1}} - {S_{\bar m,p + l + 1}} - {S_{\bar m,{\overline {p + l + 1}}}}} \right)\nonumber \\
  &\quad+ {\left( { - 1} \right)^m}\ln 2\left( {{S_{\bar lm,p + 1}} + {S_{\bar l\bar m,p + 1}} - {S_{\bar l,m + p + 1}} - {S_{\bar l,{\overline {m + p + 1}}}}} \right) \nonumber\\
  &\quad+ {\left( { - 1} \right)^{p + l}}{S_{\bar 1\bar m,{\overline {p + l + 1}}}} + {\left( { - 1} \right)^m}{S_{\bar 1\bar l,{\overline {p + m + 1}}}}.
\end{align}
\end{thm}
\pf Similarly as in the proof of formula (\ref{3.5}), we consider the Tornheim type series
\[T\left( {\bar l,\bar m;p} \right): = \sum\limits_{k,n = 1}^\infty  {\frac{{\bar H_k^{\left( l \right)}\bar H_n^{\left( m \right)}}}{{{k^p}n\left( {n + k} \right)}}} ,\quad \left( {l,m,p \in \N} \right).\]
Then with the help of formula (\ref{5.1}) we may easily deduce the result.\hfill$\square$

Taking $l=m=1,p=2r$ in (\ref{5.6}) and noting that
\[\sum\limits_{i = 1}^n {\frac{{{{\bar H}_i}}}{i}{{\left( { - 1} \right)}^{i - 1}}}  = \frac{{\bar H_n^2 + H_n^{\left( 2 \right)}}}{2},\]
then (\ref{5.6}) is equal to
\begin{align}\label{5.7}
 {S_{{{\bar 1}^3},2r + 1}} + {S_{\bar 12,2r + 1}} &= 2\bar \zeta \left( 2 \right){S_{\bar 1,2r + 1}} + 2\ln 2\left( {{S_{1\bar 1,2r + 1}} + {S_{{{\bar 1}^2},2r + 1}}} \right)\nonumber \\
  &\quad- 2\ln 2\left( {{S_{\bar 1,2r + 2}} + {S_{\bar 1,{\overline {2r + 2}}}}} \right) + 2{S_{{{\bar 1}^2},{\overline {2r + 2}}}}\nonumber \\
  &\quad- 2\sum\limits_{i = 1}^r {{{\left( { - 1} \right)}^{i - 1}}} {S_{\bar 1,i + 1}} \cdot {S_{\bar 1,2r + 1 - i}} +(-1)^{r-1} S_{\bar 1,{\overline {r + 1}}}^2,
\end{align}
Thus, we obtain the conclusion that the alternating cubic Euler sums ${S_{{{\bar 1}^3},2r + 1}}$ are reducible to alternating quadratic and linear Euler sums.

Next, we give a simple case. From \cite{BBG1994,XC2016,Xyz2016}, we know that
\begin{align}\label{5.8}
 {S_{\bar 12,3}} &= \sum\limits_{n = 1}^\infty  {\frac{{{{\bar H}_n}H_n^{\left( 2 \right)}}}{{{n^3}}}} \nonumber \\
  &= \frac{{29}}{8}\zeta \left( 2 \right)\zeta \left( 3 \right)\ln 2 - \frac{{93}}{{32}}\zeta \left( 5 \right)\ln 2 - \frac{{1855}}{{128}}\zeta \left( 6 \right) + \frac{{17}}{{16}}{\zeta ^2}\left( 3 \right)\nonumber \\
  &\quad- {S_{\bar 1,\bar 5}} + {S_{\bar 2,4}} + 4{S_{2,\bar 4}} + 8{S_{1,\bar 5}} ,
\end{align}
\begin{align}\label{5.9}
 {S_{{{\bar 1}^2},\bar 4}} &= \sum\limits_{n = 1}^\infty  {\frac{{\bar H_n^2}}{{{n^4}}}{{\left( { - 1} \right)}^{n - 1}}}\nonumber  \\
 & = \frac{{15}}{4}{\ln ^2}2\zeta \left( 4 \right) + \frac{9}{4}\zeta \left( 2 \right)\zeta \left( 3 \right)\ln 2 - \frac{{93}}{{16}}\zeta \left( 5 \right)\ln 2 + \frac{{35}}{{64}}\zeta \left( 6 \right)\nonumber \\
  &\quad- \frac{{15}}{{16}}{\zeta ^2}\left( 3 \right) + {S_{2,\bar 4}},
\end{align}
\begin{align}\label{5.10}
 {S_{{{\bar 1}^2},3}} &= \sum\limits_{n = 1}^\infty  {\frac{{\bar H_n^2}}{{{n^3}}}} \nonumber  \\
  &= 4{\rm{L}}{{\rm{i}}_5}\left( {\frac{1}{2}} \right) - \frac{1}{{30}}{\ln ^5}2 + \frac{7}{4}\zeta \left( 3 \right){\ln ^2}2 - \frac{{167}}{{32}}\zeta \left( 5 \right) + \frac{1}{3}\zeta \left( 2 \right){\ln ^3}2 \nonumber \\
  &\quad+ \frac{3}{4}\zeta \left( 2 \right)\zeta \left( 3 \right) + \frac{{19}}{8}\zeta \left( 4 \right)\ln 2,
\end{align}
\begin{align}\label{5.11}
 {S_{1\bar 1,3}} &= \sum\limits_{n = 1}^\infty  {\frac{{{H_n}{{\bar H}_n}}}{{{n^3}}}}\nonumber  \\
  &= 2{\rm{L}}{{\rm{i}}_5}\left( {\frac{1}{2}} \right) - \frac{1}{{60}}{\ln ^5}2 - \frac{{193}}{{64}}\zeta \left( 5 \right) - \frac{7}{8}\zeta \left( 3 \right){\ln ^2}2 + \frac{1}{6}\zeta \left( 2 \right){\ln ^3}2\nonumber \\
 &\quad + 4\zeta \left( 4 \right)\ln 2 + \frac{3}{8}\zeta \left( 2 \right)\zeta \left( 3 \right).
\end{align}
Setting $r=1$ in above equation (\ref{5.7}) and combining (\ref{5.8})-(\ref{5.11}), we can get the cubic sum
\begin{align}\label{5.12}
 {S_{{{\bar 1}^3},3}} &= \sum\limits_{n = 1}^\infty  {\frac{{\bar H_n^3}}{{{n^3}}}} \nonumber \\
  &= \frac{{1925}}{{128}}\zeta \left( 6 \right) - 3{\zeta ^2}\left( 3 \right) + 12{\rm{L}}{{\rm{i}}_5}\left( {\frac{1}{2}} \right)\ln 2 - \frac{{155}}{8}\zeta \left( 5 \right)\ln 2 + \frac{{27}}{8}\zeta \left( 2 \right)\zeta \left( 3 \right)\ln 2\nonumber \\
  &\quad+ \frac{{57}}{8}\zeta \left( 4 \right)\ln 2 + \frac{7}{4}\zeta \left( 3 \right){\ln ^3}2 + \zeta \left( 2 \right){\ln ^4}2 - \frac{1}{{10}}{\ln ^6}2\nonumber \\
 &\quad + {S_{\bar 1,\bar 5}} - {S_{\bar 2,4}} - 2{S_{2,\bar 4}} - 8{S_{1,\bar 5}}.
\end{align}
From \cite{XC2017}, we can find the following relation
\begin{align}\label{5.13}
{\zeta ^ \star }\left( {\bar 4,2} \right) + 4{\zeta ^ \star }\left( {\bar 5,1} \right) &= \zeta \left( {\bar 4,2} \right) + 4\zeta \left( {\bar 5,1} \right) + 5\zeta \left( {\bar 6} \right)\nonumber\\
& =  - \frac{{1105}}
{{192}}\zeta \left( 6 \right) + \frac{3}
{4}{\zeta ^2}\left( 3 \right).
\end{align}
By the definitions of alternating linear Euler sums and alternating double zeta values, we have
\[{S_{\bar p,\bar q}} = {\zeta ^ \star }\left( {\bar q,\bar p} \right),{S_{p,\bar q}} =  - {\zeta ^ \star }\left( {\bar q,p} \right),{S_{\bar p,q}} =  - {\zeta ^ \star }\left( {q,\bar p} \right).\]
Hence, the formula (\ref{5.12}) can be rewritten as
\begin{align}\label{5.14}
 {S_{{{\bar 1}^3},3}} &= \frac{{13855}}{{384}}\zeta \left( 6 \right) - \frac 3{2}{\zeta ^2}\left( 3 \right) + 12{\rm{L}}{{\rm{i}}_5}\left( {\frac{1}{2}} \right)\ln 2 - \frac{{155}}{8}\zeta \left( 5 \right)\ln 2 + \frac{{27}}{8}\zeta \left( 2 \right)\zeta \left( 3 \right)\ln 2\nonumber \\
  &\quad+ \frac{{57}}{8}\zeta \left( 4 \right)\ln 2 + \frac{7}{4}\zeta \left( 3 \right){\ln ^3}2 + \zeta \left( 2 \right){\ln ^4}2 - \frac{1}{{10}}{\ln ^6}2 + \zeta^\star\left(\bar 5,\bar 1\right)+\zeta^\star\left(4,\bar 2\right).
\end{align}
{\bf Acknowledgments.} We thank the anonymous referee for suggestions which led to improvements in the exposition.


{\small
\begin{thebibliography}{99}

\bibitem{A1997}
V. Adamchik. {\sl On Stirling numbers and Euler sums}. Journal of Computational and Applied Mathematics., 1997, {\bf 79}(1): 119-130.

\bibitem{A2000}
G.E. Andrews, R. Askey and R. Roy. {\sl Special Functions}.
Cambridge University Press., 2000: 481-532.

\bibitem{BBG1994}
D.H. Bailey, J.M. Borwein and R. Girgensohn. {\sl Experimental evaluation of Euler sums}.
Experimental Mathematics., 1994, {\bf 3}(1): 17-30.

\bibitem{B1985} B.C. Berndt. {\sl Ramanujan¡¯s Notebooks, Part I}. Springer-Verlag, New York., 1985.

\bibitem{BBG1995}
D. Borwein, J.M. Borwein and R. Girgensohn. {\sl Explicit evaluation of
Euler sums}. Proc. Edinburgh Math., 1995, {\bf 38}: 277-294.

\bibitem{BBGP1996}
J. Borwein, P. Borwein, R. Girgensohn and S.Parnes. {\sl Making sense of experimental mathematics}. Mathematical Intelligencer., 1996, {\bf 18}(4): 12-18.

\bibitem{BBBL1997}
J.M. Borwein, D.M. Bradley and D.J. Broadhurst.
{\sl Evaluations of k-fold Euler/Zagier sums:
a compendium of results for arbitrary k}. Electronic Journal of Combinatorics., 1997, {\bf 4}(2): 1-21.

\bibitem{BBBL1999}
J.M. Borwein, D.M. Bradley and D.J. Broadhurst, P. Lison${\rm \breve{e}}$k. {\sl Special values of multiple polylogarithms}. Transactions of the American Mathematical Society., 1999, {\bf 353}(3): 907-941.

\bibitem{BG1996}
J.M. Borwein and R. Girgensohn. {\sl Evaluation of triple Euler sums}. Electron. J. Combin., 1996: 2-7.

\bibitem{BZB2008}
J.M. Borwein, I.J. Zucker and J. Boersma. {\sl The evaluation of character Euler double sums}. Ramanujan J., 2008, {\bf 15} (3): 377-405.

\bibitem{C1974}
L. Comtet. {\sl Advanced combinatorics}. Boston: D Reidel Publishing Company., 1974.

\bibitem{CB1994}
R.E. Crandall and J.P. Buhler. {\sl On the evaluation of Euler sums}. Experimental Mathematics., 1994, {\bf 4}(3): 275-285.

\bibitem{ELY2005}
M. Eie, W. Liaw and F. Yang. {\sl On evaluation of generalized Euler sums of even weight}. Int. J. Number Theory., 2005, {\bf 225}(1): 225-242.

\bibitem{EW2012}
 M. Eie and C. Wei. {\sl Evaluations of some quadruple Euler sums of even weight}. Functions et Approximatio, 2012, {\bf 46}(1): 63--67.

\bibitem{FS1998}
P. Flajolet and B. Salvy. {\sl Euler sums and contour integral representations}. Experimental Mathematics., 1998, {\bf 7}(1): 15--35.

\bibitem{F2005}
P. Freitas. {\sl Integrals of polylogarithmic functions, recurrence relations, and associated Euler sums}. Mathematics of Computation., 2005, {\bf 74}(251): 1425-1440.

\bibitem{H1992}
M.E. Hoffman. {\sl Multiple harmonic series}. Pacific J. Math., 1992, (152): 275-290.

\bibitem{M2014}
I. Mez${\rm \ddot{o}}$. {\sl Nonlinear Euler sums}. Pacific J. Math., 2014, {\bf 272}: 201-226.

\bibitem{KP2013}
Kh. Hessami Pilehrood, T. Hessami Pilehrood and R. Tauraso. {\sl New properties of multiple harmonic sums modulo $p$ and $p$-analogues of Leshchiner's series}. Trans. Amer. Math. Soc., 2014, {\bf 366}(6): 3131-3159.

\bibitem{O2008}
Y.L. Ong, M. Eie and W.C. Liaw. {\sl Explict evaluation of Triple Euler sums}. Int. J. Number Theory., 2008, {\bf 4}(3): 437-451.

\bibitem{O2010}
Y.L. Ong, M. Eie and C. Wei. {\sl Explicit evaluations of quadruple Euler sums}. Acta Arithmetica., 2010, {\bf 144}: 213--230.

\bibitem{M2000}
H.N. Minh and M. Petitot. {\sl Lyndon words, polylogarithms and the Riemann $\zeta$ function}. Discrete Mathematics., 2000, {\bf 217}: 273-292.

\bibitem{SX2017}
X. Si and C. Xu. {\sl Evaluations of some quadratic Euler sums}. arXiv:1701.00389 [math.NT].

\bibitem{T2004}
H. Tsumura. {\sl Combinatorial relations for Euler-Zagier sums}. Acta. Arith., 2004, {\bf 111}: 27--42.

\bibitem{W2017}
W. Wang and Y. Lyu. {\sl Euler sums and Stirling sums}. J. Number Theory., in press.

\bibitem{XC2016}
C. Xu and Y. Cai. {\sl On harmonic numbers and nonlinear Euler sums}. arXiv:1609.04924 [math.NT].

\bibitem{Xu2016}
C. Xu and J. Cheng. {\sl Some results on Euler sums}. Functions et Approximatio., 2016, {\bf 54}(1): 25-37.

\bibitem{Xu2017}
C. Xu. {\sl Multiple zeta values and Euler sums}. J. Number Theory., 2017, {\bf 177}: 443-478.

\bibitem{X2017}
C. Xu and Z. Li. {\sl Tornheim type series and nonlinear Euler sums}. J. Number Theory., 2017, {\bf 174}: 40-67.

\bibitem{X2016}
C. Xu, Y. Yan and Z. Shi. {\sl Euler sums and integrals of polylogarithm functions}. J. Number Theory., 2016, {\bf 165}: 84-108.

\bibitem{Xyz2016}
C. Xu, Y. Yang and J. Zhang. {\sl Explicit evaluation of quadratic Euler sums}. Int. J. Number Theory., 2017, {\bf 13}(3): 655-672.

\bibitem{XMZ2016}
C. Xu, M. Zhang and W. Zhu. {\sl Some evaluation of harmonic number sums}. Integral Transforms and Special Functions., 2016, {\bf 27}(12): 937-955.

\bibitem{XC2017}
C. Xu. {\sl Evaluations of nonlinear Euler sums of weight ten}. arXiv:1703.00254v2[math.NT].

\bibitem{DZ1994}
D. Zagier.{\sl Values of zeta functions and their applications}. First European Congress
of Mathematics, Volume II, Birkhauser, Boston., 1994, (120): 497-512.

\bibitem{DZ2012}
D. Zagier. {\sl Evaluation of the multiple zeta values $\zeta(2,...,2,3,2,...,2)$}. Annals of Mathematics., 2012, {\bf 2}(2): 977-1000.

\bibitem{Z2015}
J. Zhao. {\sl Multiple Zeta functions, multiple polylogarithms and their special values}. World Scientific, 2015.

\end{thebibliography}}
 \end{document}